\documentclass{amsart}

\usepackage{amsmath,amssymb,amsthm}

\newtheorem{theorem}{Theorem}
\newcommand{\bt}{\begin{theorem}}
\newcommand{\et}{\end{theorem}}
\newtheorem{lemma}{Lemma}
\newcommand{\bl}{\begin{lemma}}
\newcommand{\el}{\end{lemma}}
\newtheorem{corollary}{Corollary}
\newcommand{\bc}{\begin{corollary}}
\newcommand{\ec}{\end{corollary}}
\newtheorem{problem}{Problem}
\newcommand{\bprob}{\begin{problem}}
\newcommand{\eprob}{\end{problem}}
\newcommand{\beq}{\begin{equation}}
\newcommand{\eeq}{\end{equation}}
\newcommand{\benum}{\begin{enumerate}}
\newcommand{\eenum}{\end{enumerate}}
\newcommand{\N}{\ensuremath{ \mathbf N }}
\newcommand{\Z}{\ensuremath{\mathbf Z}}
\newcommand{\mcp}{\ensuremath{ \mathcal P}}
\newcommand{\mcq}{\ensuremath{ \mathcal Q}}
\newcommand{\diam}{\ensuremath{\text{diam}}}
\newtheorem*{NoNumPrachar}{Prachar's Theorem}

\title[Arithmetic diameter of the integers]{Geometric group theory and arithmetic diameter}

\author{Melvyn B.  Nathanson}
\address{Lehman College (CUNY), 
Bronx, NY 10468, and 
CUNY Graduate Center,
New York, NY 10016}
\email{melvyn.nathanson@lehman.cuny.edu}

\dedicatory{To K. Gy\H  orgy and A. S\' ark\" ozy on their 70th birthdays\\ and A. Peth\H o and J. Pintz on their 60th birthdays}


\date{\today}

\begin{document}

\maketitle

\begin{abstract}
Let $X$ be a group with identity $e$, let $A$ be an infinite set of generators for $X$, and let $(X,d_A)$ be the metric space with the word metric $d_A$ induced by $A$.  If the diameter of the space is infinite, then for every positive integer $h$ there are infinitely many elements $x\in  X$ with $d_A(e,x)=h$.  
It is proved that if \mcp\ is a nonempty finite set of prime numbers and  $A$ is the set of positive integers whose prime factors all belong to \mcp, then the metric space $(\Z,d_A)$ has infinite diameter.  
Let $\lambda_A(h)$ denote the smallest positive integer $x$  with $d_A(e,x)=h$. 
It is an open problem to compute $\lambda_A(h)$ and estimate its growth rate.
\end{abstract}

\section{Word metrics on infinitely generated groups}
There is a nice interaction between geometric group theory and additive number theory, and concepts from each have fruitful analogs in the other. 
Indeed, a large part of geometric group theory can be considered as a kind of ``nonabelian additive number theory.''
We shall consider  growth in groups with a fixed infinite set of generators. 
The group of integers, of course, is of unique interest to number theorists, and we describe some problems and results about additive bases of infinite order that arise in  geometry.

Let $\N_0$ denote the nonnegative integers and \Z\ the integers.

Let $X$ be a group, written multiplicatively, with identity $e$. 
For $A \subseteq X$, 
we define $A^{-1} = \left\{ a^{-1} : a\in A\right\}$.
A subset $A$ of $X$ is a  \emph{generating set} for $X$ if every element of $X$ can be written as a finite product of elements of $A$ and their inverses.  
Equivalently, $A$ generates $X$ as a group if and only if $A\cup A^{-1}$ generates $X$ as a semigroup.
The \emph{length} of an element $x \in X$ with respect to $A$ is
\begin{align*}
\ell_A(x) & = \min\{k : x =  a_1^{\varepsilon_1} a_2^{\varepsilon_2} \cdots a_k^{\varepsilon_k} : a_i \in A \text{ and } \varepsilon_i \in \{1,-1\} \text{ for } i=1,2,\ldots, k\}  \\
&  = \min\{k : x =  a_1a_2\cdots a_k : a_i \in A  \cup A^{-1}
\text{ for } i=1,2,\ldots, k\}.
\end{align*}
Thus, in a group, $\ell_A(x)$ is the number of letters in the shortest word in the alphabet $A \cup A^{-1}$ that spells $x$.  
The length function of a group  has the following properties:
\begin{enumerate}
\item
$\ell_A(x) \geq 0$ for all $x\in X$
\item
$\ell_A(x) = 0$ if and only if $x=e$
\item
For all $x,y \in X$,
\[
\ell_A(xy) \leq \ell_A(x)+\ell_A(y)
\]
\item
For all $x\in X$,
\[
\ell_A\left(x^{-1} \right) = \ell_A(x).
\]
\end{enumerate}
These properties of the length function imply that the function $d_A:X \times X\rightarrow \N_0$ defined by 
\[
d_A(y,x) = \ell_A\left(y^{-1}x\right)
\]
is a distance function on the group $X$.  
We call $d_A$ the  \emph{word metric} induced by $A$.  
Geometric properties of the metric space $(X,d_A)$ 
are often related to combinatorial and additive number theoretic properties of the set $A$.  

The definition of the length function also implies that, for all $x\in X$ and $a \in A \cup A^{-1}$,  
\beq    \label{AD:LengthIneq-1}
\ell_A(x)-1 \leq \ell_A(ax) \leq \ell_A(x) + 1
\eeq
and, if $\ell_A(x) = h$ and $x=a_1a_2 \cdots a_h$ with $a_i \in A \cup A^{-1}$ for $i=1,2,\ldots, h$, then 
\beq    \label{AD:LengthIneq-2}
\ell_A(a_{j+1}a_{j+2}\cdots a_k) = k-j
\eeq
for all $0 \leq j < k \leq h$.

For every nonnegative integer $h$, 
the \emph{sphere} of radius $h$ with respect to $A$ is the set
\beq \label{AD:defineSphere}
S_A(h) = \left\{ x\in X_0:\ell_A(x) = h \right\}
\eeq
and the  \emph{closed ball} of radius $h$ with respect to $A$ is the set
\beq \label{AD:defineBall}
\widetilde{B_A}(h) = \left\{ x\in X_0:\ell_A(x) \leq h \right\}
= \bigcup_{h_1=0}^h S_A(h_1).
\eeq

Most research in geometric group theory concerns finitely generated groups.  
In number theory, however, it is often important to study a group $X$ with an infinite set $A$ of generators.    
The generating set $A$ is not necessarily minimal, and it is irrelevant whether or not the group can be finitely generated.  

The following theorem is fundamental for the geometry of groups with an infinite number of generators.  
Recall that the \emph{diameter} of a metric space $(X,d)$ is 
\[
\diam(X) = \sup\{ d(x,y):x,y\in X\}.
\]

\bt  \label{AD:theorem:spheres}
Let $X$ be a group and let $A$ be an infinite generating set for $X$.  If the sphere $S_A(h)$ is finite for some positive integer $h$, then $S_A(k) = \emptyset$ for all $k  > h$.  
Equivalently, if $\diam(X) = \infty$, then $S_A(h)$ is an infinite set for 
every positive integer $h$.  
\et

\begin{proof}
Let $\ell_A$ be the length function associated with the generating set $A$.  
Let $h$ be a  positive integer such that $S_A(h)$ is finite.
The sphere $S_A(1) = \left( A \cup A^{-1}\right)\setminus \{ e\}$ is infinite, and so $h \geq 2.$  
If $S_A(k) \neq \emptyset$  for some integer $k  > h$, then there exist elements $a_1,\ldots, a_k \in A \cup A^{-1}$ such that 
\[
\ell_A(a_1a_2 \cdots a_k)= k.
\]
We  define a function $f:A \rightarrow \{1,\ldots, k \}$ as follows.
By~\eqref{AD:LengthIneq-1}, for every $a\in A$ we have
\[
\ell_A(aa_1a_2 \cdots, a_k)  \geq k -1 \geq h.
\]
It follows that
there is at least one integer $t \in \{1,\ldots, k \}$ such that 
\[
\ell_A(aa_1a_2 \cdots a_t) = h.
\]
Choose any  $t$ with $\ell_A(aa_1a_2 \cdots a_t) = h$ 
and let $f(a) = t$.  
Because the set $A$ is infinite and 
$A = \bigcup_{i=1}^k f^{-1}(i)$, 
it follows from the pigeonhole principle that $f^{-1}(t)$ is an infinite set for some 
$t \in  \{1,\ldots, k \}$, and so there exist infinitely many elements $a\in A$ such that $aa_1a_2 \cdots a_t$ has length $h$.  
If $a\neq a'$, then $aa_1a_2 \cdots a_t \neq a'a_1a_2 \cdots a_t$
and  the finite set $S_A(h)$ must contain the infinite subset 
$\{ aa_1a_2 \cdots a_t  : a \in f^{-1}(t)\}$, which is absurd.  
Thus, if $S_A(h)$ is finite for some $h\geq 2$, then $S_A(k)=\emptyset$ for all $k > h$.  
This completes the proof.  
\end{proof}

\section{Sum and difference sets for additive groups}

If $X$ is an additive abelian group with identity 0, then for every subset $A$ of $X$ and for every positive integer $h$, 
we define the $h$-fold \emph{sum and difference set}
\[
h^{\pm}A = \left\{ \sum_{i=1}^h \varepsilon_i a_i : 
\varepsilon_i \in \{ 1,-1\} \text{ and } a_i \in A 
\text{ for } i=1,\ldots, h\right\}.
\]
We define $0^{\pm}A = \{0\}$. 
If $A$ generates the group $X$, then, in geometric group theory, the set 
\[
\widetilde{B_A}(h) = \bigcup_{h_1 \leq h} h_1^{\pm}A = h^{\pm}(A \cup \{ 0\})
\]  
is the closed ball of radius $h$ in  the metric space $(X,d_A)$ 
with respect to the word metric $d_A$ induced by the generating set $A$.  
The generating set $A$ is called a \emph{basis of order $h$} for $X$  
if $\widetilde{B_A}(h) = X$.
The set $A$ is a \emph{basis} for $X$ if 
$\widetilde{B_A}(h) = X$ for some positive integer $h$.

The sum and difference sets $h^{\pm}A$ are finite for all positive integers $h$ if and only if the set $A$ is finite.  
Moreover, if $A$ is finite, 
then 
\beq    \label{AD:SumsetIneq-1}
|h^{\pm}A| \leq {2|A| +h-1 \choose h} \leq (2|A|)^h
\eeq
and
\beq    \label{AD:SumsetIneq-2}
|\widetilde{B_A}(h)| \leq {2|A| +h \choose h} \leq (2|A|+1)^h
\eeq
where $|A|$ denotes the cardinality of the set $A$.

\bt    \label{AD:theorem:DiameterBasis}
Let $A$ be a generating set for the group $X$.  
The set $A$ is a basis for $X$ if and only if the diameter of the metric space $(X,d_A)$ is finite. 
\et

\begin{proof}
The diameter of the metric space $(X,d_A)$ is finite if and only if $\widetilde{B_A}(h) = X$ for some positive integer $h$.
\end{proof}

\section{Arithmetic diameter}
A subset $A$ of the additive abelian group \Z\ is a generating set for \Z\ if and only if $\gcd(A) = 1$.  
The \emph{arithmetic diameter} of a set $A$ of relatively prime integers is the diameter of the metric space $(\Z,d_A)$.
It is, in general, a difficult problem to compute or to estimate the arithmetic diameter of a generating set, or even to describe the  infinite generating sets whose arithmetic diameters are  finite.  
In this section we prove that if \mcp\ is a nonempty finite set of prime numbers and if $A$ is the set of positive integers all of whose prime factors belong to $A$, then $\diam(\Z,d_A) = \infty$.
We need the following beautiful result from analytic number theory.  

\begin{NoNumPrachar}
Let $\delta(n)$ denote the number of prime numbers $p$ such that $p-1$ divides $n$.  For every positive real number $c < (\log 2)/2$ and for every positive integer $k_0$, there exist infinitely many positive integers $n$ such that 
\[
\delta(n) > \exp\left( \frac{c\log n}{\log\log n}\right) + k_0.
\]
\end{NoNumPrachar}

\begin{proof}
See Prachar~\cite{prac55}, or Adleman, Pomerance, and Rumely~\cite[Prop. 10]{adle-pome-rume83}.
\end{proof}

\bt    \label{AD:theorem:smooth}
Let $\mcp$ be a nonempty finite set of prime numbers and let $A$ be the set of positive integers that are products of primes in \mcp, that is,
\[
A = \left\{ \prod_{p \in \mcp} p^{v_p} : v_p \in \N_0 \right\}.
\]
Let $d_A$ be the word metric induced by $A$ on the group \Z.  
The diameter of the metric space $(\Z, d_A)$ is infinite.  
Equivalently, $h^{\pm}(A \cup \{ 0\}) \neq \Z$ for every positive integer $h$.
\et

\begin{proof}
Let   \mcq\ be a finite set of prime numbers that is disjoint from \mcp, 
and let 
\[
Q = \prod_{q\in \mcq} q.
\]
If $a \in A$ and $q \in \mcq$, then $\gcd(a,q) = 1$ and  
\[
a^{q-1}  \equiv 1 \pmod{q}.
\]
Let $n$ is an integer that is divisible by $q-1$ for all $q \in \mcq$.  
Then
\[
a^{n}  \equiv 1 \pmod{q}
\]
for all $q\in \mcq$, and so
\beq  \label{AD:Qcong} 
a^{n}  \equiv 1 \pmod{Q}.
\eeq

Let $\Z/Q\Z$ be the ring of congruence classes modulo $Q$.
For every integer $r$, we denote the congruence class of $r$ modulo $Q$ by $\overline{r} = r +Q\Z$.  
For every set $R$ of integers, we define  
\[
\overline{R} = \left\{ \overline{r} : r \in R  \right\} \subseteq \Z/Q\Z.
\]
Consider the set of integers
\[
R_Q = \left\{\prod_{p\in \mcp} p^{j_p} : j_p \in \{0,1,\ldots, n-1\}  \right\}.
\]
We have $R_Q \subseteq A$ and so $\overline{R_Q} \subseteq \overline{A}$.

Let $a = \prod_{p\in \mcp} p^{v_p} \in A$.  
By the division algorithm, for each prime $p \in \mcp$ there are integers $w_p$ and $j_p$ with $j_p \in \{0,1,\ldots, n-1\}$ such that 
\[
v_p = w_p n + j_p.
\]
Defining
\[
a_1 = \prod_{p\in \mcp} p^{w_p} \in A
\]
and
\[
r =  \prod_{p\in \mcp} p^{j_p} \in R_Q
\]
we obtain 
\[
a = \prod_{p\in \mcp} p^{v_p}  = \prod_{p\in \mcp} p^{w_p n + j_p} 
= \left(  \prod_{p\in \mcp} p^{w_p } \right)^n  \prod_{p\in \mcp} p^{j_p}
= a_1^n r.
\]
Congruence~\eqref{AD:Qcong} implies that  
\[
a \equiv r \pmod{Q}
\]
and so $\overline{a} = \overline{r} \in \overline{R_Q}$ for all $a\in A$, that is, 
$\overline{A} \subseteq \overline{R_Q}$, hence  
$\overline{A} = \overline{R_Q}$.
Let
\[
|\mcp| = k_0.
\]
Then 
\[
|\overline{A}| = |\overline{R_Q}| \leq |R_Q| = n^{k_0}
\]
and, by~\eqref{AD:SumsetIneq-2},
\[
\left| h^{\pm}\left( \overline{ R_Q} \cup\{ \overline{0} \} \right) \right| 
\leq \left( 2 |\overline{R_Q}|+1 \right)^h 
\leq \left( 2 n^{k_0}  +1 \right)^h 
<  c_0 n^{k_0h}
\]
for $c_0 = 3^h$.  
If $\diam(\Z,d_A) = h$, then 
\[
h^{\pm}(A\cup \{0\}) = \Z
\]
and so
\[
 \Z/Q\Z = 
 \overline{ h^{\pm}(A\cup \{ 0 \} ) } 
 = h^{\pm} \left(    \overline{A} \cup \{ \overline{0} \}    \right) 
= h^{\pm}\left( \overline{R_Q} \cup\{ \overline{0} \}  \right) 
\]
and
\[
Q = | \Z/Q\Z| = \left| h^{\pm}\left( \overline{R_Q} \cup\{ \overline{0} \}\right) \right| < c_0n^{k_0h}.
\]
Conversely, if we can construct a finite set $\mcq$ of prime numbers that is disjoint from \mcp\ and an integer $n$ that is divisible by $q-1$ for every $q\in \mcq$ such that  $c_0 n^{k_0h} < Q$, then $\diam(\Z,d_A) \neq h$.
This is what we shall do.

Applying Prachar's theorem with $k_0 = |\mcp|$, we obtain arbitrarily large integers $n$ with the property that there are at least
\[
\delta(n) - k_0 > \exp\left( \frac{c\log n}{\log\log n}\right) 
\]
prime numbers $q$ such that $q-1$ divides $n$ 
and $q\notin \mcp$. 
For each such $n$, let $\mcq_n$ be the set of these primes and let $Q_n$ be the product of the primes in $\mcq_n$.   Then 
\begin{align*}
Q_n & \geq 2^{ \delta(n) - k_0  } \\
& > \exp{\left( \log 2 \exp\left( \frac{c\log n}{\log\log n}\right) \right)} \\
& >  c_0 n^{k_0h}
\end{align*}
for $n$ sufficiently large.
This completes the proof.
\end{proof}

\bt     \label{AD:theorem:BigSphere}
Let $\mcp$ be a nonempty finite set of prime numbers and let $A$ be the set of positive integers that are products of primes in \mcp.  
Let $h$ be a positive integer 
and let $S_A(h) = \{n \in \Z : \ell_A(n) = h\}.$ 
Then $S_A(h)$ contains infinitely many integers for every positive integer $h$.
\et

\begin{proof}
This follows immediately from Theorem~\ref{AD:theorem:spheres}.
\end{proof}

\bl     \label{AD:lemma:subsetAD}
Let $A$ and $A'$ be generating sets for \Z, with induced word metrics $d_A$ and $d_{A'}$, respectively.  
If $A \subseteq A'$, then $d_{A}(x,y) \geq d_{A'}(x,y)$ 
for all $x,y \in \Z$.  
In particular, if the arithmetic diameter of the set $A'$ is infinite, then the  arithmetic diameter of the set $A$ is infinite.
\el

\begin{proof}
This follows immediately from the definitions of word metric and  diameter.  
\end{proof}

\bt    \label{AD:theorem:Exponentials}
Let $\{a_1,\ldots, a_k\}$ be a set of positive integers.  
The arithmetic diameter of the set  
$A = \bigcup_{i=1}^k \left\{a_i^j : j \in \N_0  \right\}$ is infinite. 
Moreover, $\{ n \in \Z: \ell_A(n) = h \}$ is infinite for every positive integer $h$.
\et

\begin{proof}
Let \mcp\ be the finite set consisting of all prime numbers $p$ such that that $p$ divides $a_i$ for some $i=1,\ldots, k$.  
Let $A'$ be the set of positive integers all of whose prime factors belong to \mcp.   Then $A \subseteq A'$.
The set $A'$ has infinite arithmetic diameter by Theorem~\ref{AD:theorem:smooth}, and so $A$ has infinite arithmetic diameter by Lemma~\ref{AD:lemma:subsetAD}.
The sphere $S_A(h) = \{ n \in \Z: \ell_A(n) = h \}$ is infinite for every positive integer $h$ by Theorem~\ref{AD:theorem:BigSphere}.  
This completes the proof.
\end{proof}

\section{Open problems}

The following result is a  special case of Theorem~\ref{AD:theorem:Exponentials}.

\bt   \label{AD:theorem:Exponentials-23}
For every positive integer $h$ there are infinitely many positive integers $N$ such that, for some nonnegative integers $j$ and $k$ with $h=j+k$,  the exponential diophantine equation 
\[
N = \pm 2^{v_1} \pm 2^{v_2} \pm \cdots \pm 2^{v_j} \pm 3^{w_1} \pm 3^{w_2} \cdots \pm 3^{w_k}
\]
has a solution in nonnegative integers $v_1,\ldots, v_j, w_1,\ldots, w_k$, but the corresponding equation with $h$ replaced by any positive integer $h' < h$ has no solution.
\et

\begin{proof}
Apply Theorem~\ref{AD:theorem:Exponentials} to the set $\{2,3\}$.
\end{proof}

\bprob
Theorem~\ref{AD:theorem:Exponentials-23} implies that for every positive integer $h$ there is a smallest number $N = \lambda_{2,3}(h)$ that can be written as the sum or difference of exactly $h$ but no fewer powers of 2 and 3.
For example, $ \lambda_{2,3}(1)=1$, $ \lambda_{2,3}(2)=5$, 
$ \lambda_{2,3}(3) \geq 150$, and $ \lambda_{2,3}(3) = 150$ 
if and only if at least one of the exponential diophantine equations 
$149 = |2^x \pm 3^y|$ and $151 = |2^x \pm 3^y|$ has a solution in nonnegative integers.  
Eamonn O'Brien (personal communication) showed that $ \lambda_{2,3}(4) \geq 393$
Compute the function $ \lambda_{2,3}(h)$.  Is it increasing?  Estimate its rate of growth.  
\eprob

\bprob
Let \mcp\ be a nonempty finite set of prime numbers, and let $A$ be the set of positive integers that are products of powers of primes in \mcp.   By Theorem~\ref{AD:theorem:smooth},
for every positive integer $h$ there are infinitely many  positive integers of length $h$.  Let $\lambda_{\mcp}(h)$ denote the smallest  positive integer of length $h$.  Determine the properties of this function.  Estimate its growth rate.
\eprob

\bprob
Let $(a_i)_{i=1}^{\infty}$ be a strictly increasing sequence of positive integers, and let
\[
A = \bigcup_{i=1}^{\infty} \{ \pm a_i^j : j=i,i+1.i+2,\ldots \}.
\]
Is $A$ a basis of finite order for \Z?
The condition $j \geq i$ is used to avoid Waring's problem.
\eprob

\bprob
Let $a$ and $b$ be integers with $a\geq 2$ and $b\geq 2$, and let $(\Z,d_a)$ and $(\Z,d_b)$ be the metric spaces with the word metrics induced by the generating sets $\{a^i:i\in \N_0\}$ and 
$\{ b^i : i\in \N_0\}$, respectively.  
Nathanson~\cite{nath11b} proved that the word metrics $d_a$ and $d_b$ are bi-Lipschitz equivalent if and only if there exist positive integers $m$ and $n$ such that $a^m = b^n$.  
It is not know, however, under what conditions the metric spaces 
$(\Z,d_a)$ and $(\Z,d_b)$ are bi-Lipschitz equivalent.  It is a problem of Richard E. Schwartz in metric geometry to determine if the metric spaces $(\Z,d_a)$ and $(\Z,d_b)$ are quasi-isometric.  This is unknown even in the case $a=2$ and $b=3$.
\eprob

\bprob
Let $A$ and $B$ be infinite generating sets for \Z.  Determine necessary and sufficient conditions for the metric spaces $(\Z,d_A)$ and $(\Z,d_B)$ to be quasi-isometric.
\eprob

\section{Notes}
 In November, 2008, in a lecture in the Princeton Discrete Mathematics Seminar,  I discussed the bi-Lipschitz inequivalence of the metrics on the spaces $(\Z,d_2)$ and $(\Z,d_3)$.
 After the talk, John Conway suggested we try to compute the function $\lambda_{2,3}(h)$.  
We failed, but the effort led to the problems discussed in this paper.

Theorem~\ref{AD:theorem:spheres} appears in Nathanson~\cite{nath11c}.
I learned recently from Hannah Alpert that essentially the same result, but written in the language of graph theory, is in a paper of Wesley Pegden~\cite{pegd06}.

Theorem~\ref{AD:theorem:Exponentials-23} was first proved by Gautam Chinta (personal communication).    
The proof of Theorem~\ref{AD:theorem:smooth}  is based on Chinta's idea.  
I learned from the preprint of Hajdu and Tijdeman~\cite{hajd-tijd11}  that this result  was  proved by Jarden and Narkiewicz~\cite{jard-nark07} in the language of algebraic number theory.   

It is useful to be polylingual in mathematics.

\def\cprime{$'$} \def\cprime{$'$} \def\cprime{$'$}
\providecommand{\bysame}{\leavevmode\hbox to3em{\hrulefill}\thinspace}
\providecommand{\MR}{\relax\ifhmode\unskip\space\fi MR }
\providecommand{\MRhref}[2]{%
  \href{http://www.ams.org/mathscinet-getitem?mr=#1}{#2}
}
\providecommand{\href}[2]{#2}

\end{document}